\documentclass[matann]{svjour}
\usepackage{amsmath,amssymb,amsfonts,amscd,indentfirst}
\usepackage[dvipdfm]{hyperref}
\usepackage{graphics}

\def\cE{{\mathcal E}}
\def\cF{{\mathcal F}}

\def\cW{{\mathcal W}}

\def\bg{{\mathbf g}}

\newtheorem{prop}{Proposition}[section]
\newtheorem{theo}[prop]{Theorem}

\newtheorem{coro}[prop]{Corollary}
\newtheorem{rema}[prop]{Remark}

\newtheorem{defi}[prop]{Definition}

\def\begeq{\begin{equation}}
\def\endeq{\end{equation}}

\def\and{\quad{\rm and}\quad}

\def\<{\langle}
\def\>{\rangle}

\def\Dint{\displaystyle\int}
\def\Dfrac{\displaystyle\frac}

\begin{document}
\title{Eigenvalues and energy functionals with monotonicity formulae under
Ricci flow}
\author{Jun-Fang Li}
\institute{\textsc{Jun-Fang Li}\at Mathematical Sciences Research
Institute, Berkeley, CA 94720, USA\\
email : junfangl@msri.org}
\date{Received: date / Revised version: date}
%
\titlerunning{Monotonicity formulae under Ricci flow}
\authorrunning{Jun-Fang Li}
\maketitle
\begin{abstract}
In this note, we construct families of functionals of the type of
$\mathcal{F}$-functional and $\mathcal{W}$-functional of Perelman.
We prove that these new functionals are nondecreasing under the
Ricci flow. As applications, we give a proof of the theorem that
compact steady Ricci breathers must be Ricci-flat. Using these new
functionals, we also give a new proof of Perelman's no non-trivial
expanding breather theorem. Furthermore, we prove that compact
expanding Ricci breathers must be Einstein by a direct method. In
this note, we also extend Cao's methods of eigenvalues\cite{C} and
improve their results. \subclass{53C44; 35K55}
\end{abstract}
\section{\textnormal{\bf Introduction}}\label{Section Introduction}
Let ($M,\bg$) be a closed Riemannian manifold. In \cite{P}, Perelman
introduced a functional
\begin{equation} \label{Intro def1}
 \mathcal{F}(\mathbf{g},f) = \Dint_M
(R + |\nabla f|^2)e^{-f}d\mu.
\end{equation}
If ($M,\bg$) is a solution to the Ricci flow equation, Perelman
proved that the $\mathcal{F}$-functional is nondecreasing under the
Ricci flow and the Ricci flow can be viewed as the gradient flow of
this functional. He proved that under the following coupled system:
\begin{equation}
\left\{
\begin{array}{rll}
\frac{\partial}{\partial t}\mathbf{g}_{ij} &= -2R_{ij}\\
\frac{\partial}{\partial t}f &= -\Delta f - R + |\nabla f|^2,
\end{array}
\right.\label{sec1.1}
\end{equation}
the $\mathcal{F}$-functional is nondecreasing. More precisely,
\begin{equation}\label{sec1.2}
\frac{d}{d t}\mathcal{F} = 2\Dint_M |R_{ij} + \nabla_i\nabla_j f|^2
e^{-f}d\mu \ge 0.
\end{equation}

If one defines
\[
\lambda(\mathbf{g}) = \inf\mathcal{F}(\mathbf{g},f),
\]
where the infimum is taken over all the smooth $f$ which satisfies

\[
\Dint_M e^{-f}d\mu = 1,
\]
then $\lambda(\mathbf{g})$ is the lowest eigenvalue of the operator
\[
-4\Delta + R,
\] and the nondecreasing of the $\mathcal{F}$ functional implies the
nondecreasing of $\lambda(\mathbf{g})$. As an application, Perelman
was able to show that there is no nontrivial steady or expanding
Ricci breathers on closed manifolds.\\

In \cite{C}, X.D. Cao considered the eigenvalues of the operator
\[
-\Delta + \frac{R}{2},
\]
on manifolds with {\em nonnegative curvature operator}. They showed
that the eigenvalues of these manifolds are nondecreasing along the
Ricci flow. Using the monotonicity of the eigenvalues, they proved
that the only steady Ricci breather with {\em nonnegative curvature
operator} is the trivial one which is Ricci-flat.\\

In this note, we study monotonicity formulae of various energy
functionals and the eigenvalues of the operator $-4\Delta + kR$. In
section \ref{Section A remark on the assumption of positive Ricci
curvature operator},\footnote{The results in section \ref{Section A
remark on the assumption of positive Ricci curvature operator} are
relatively independent of the rest of the paper and can be treated
separately.} we study a monotonicity formula of eigenvalues of
$-\Delta +\frac{R}{2}$, which improves the result in \cite{C} based
on the same technique Cao has used. The following is one of the main
theorems in this paper :

\begin{theo}\label{thmljf1}
Let $\bg(t)$, $t\in[0,T)$, be a solution to the Ricci flow on a
closed Riemannian manifold $M^n$. Assume that there is a
$C^1$-family of smooth functions $f(t)>0$, which satisfy
\begin{equation}\label{thmljf1 rev equ1}
\lambda(t)f(t) = -\Delta_{\bg(t)}f(t) + \frac{1}{2}R_{\bg(t)}f(t)
\end{equation}
\begin{equation}\label{thmljf1 rev equ2}
\Dint_M f^2(t)d\mu_{\bg(t)} = 1
\end{equation}
where $\lambda(t)$ is a function of $t$ only. Then
\begin{equation} \label{steady positive equ11}
\begin{array}{rll}
2\frac{d}{dt}\lambda(t) &= 4\Dint R_{ij}\nabla^if\nabla^jf d\mu +
2\Dint|Rc|^2f^2d\mu\\
&=\Dint |R_{ij}+\nabla_i\nabla_j \varphi|^2\ e^{-\varphi}d\mu +
\Dint
|Rc|^2\ e^{-\varphi}d\mu\\
&\ge 0,
\end{array}
\end{equation}
where $\varphi$ satisfies $e^{-\varphi} = f^2$.
\end{theo}

 In section \ref{Section From first variation formula to
Entropy functionals}, we introduce a functional $\mathcal{E}$ and
derive its first variation along the Ricci flow. In section
\ref{Section Entropy functionals}, we define new functionals
$\mathcal{F}_k$ and we prove that although the Ricci flow is not a
gradient flow of $\mathcal{F}_k$ in the sense of Perelman's,
$\mathcal{F}_k$ are still nondecreasing under the Ricci flow. The
monotonicity is strict unless we are on a Ricci-flat solution.

 In section \ref{section
steady}, by using the monotonicity of the eigenvalues, we rule out
compact steady Ricci breathers.

In section \ref{Section New entropy formulae over expanders}, we
introduce a family of new functionals $\mathcal{W}_e$ and
$\cW_{ek}$. These functionals are not scale invariant in contrast to
Perelman's $\mathcal{W}$-entropy. We obtain their first variations :

\begin{theo}\label{Entropy formulae over Expanders main}
Under the following coupled system
$$\left\{
\begin{array}{rll}
\Dfrac{\partial}{\partial t}\mathbf{g}_{ij} &= -2R_{ij}\\
\frac{\partial}{\partial t}f &= -\Delta f +|\nabla f|^2- R\\
\Dfrac{d}{dt}\tau &= 1,
\end{array}
\right.
$$
we have the following monotonicity,\footnote{See the definition of
$\cW_{ek}$ in (\ref{Entropy formulae over Expanders equ1}) }
$$\begin{array}{rll}
\Dfrac{d}{d t} \cW_{ek} =& 2\tau^2\Dint_M
\big|R_{ij}+\nabla_i\nabla_j
f+\Dfrac{1}{2\tau}\mathbf{g}_{ij}\big|^2e^{-f}d\mu \\
&\quad+ 2(k-1)\tau^2\Dint_M
\big|R_{ij}+\Dfrac{1}{2\tau}\mathbf{g}_{ij}\big|^2e^{-f}d\mu,
\end{array}
$$ where $k\ge 1$. The
monotonicity is strict unless we are on an expanding gradient Ricci
soliton ($k=1$) or an Einstein manifold ($k>1$).
\end{theo}

In section \ref{Section No expanding breathers theorem}, we use the
non-scale-invariant functional $\mathcal{W}_e$ to prove Perelman's
no non-trivial expanding breather theorem.

 In section \ref{Section Expanding Ricci breathers
are Einstein}, we use $\cW_{ek}$ functionals to rule out
non-Einstein compact expanding Ricci breathers.

\begin{rema}
The no non-Einstein steady or expanding Ricci breather theorem was
first proved by T. Ivey in \cite{I}. We give a different proof here
which is based on a global functional method and does not use
maximal principles.
\end{rema}

\begin{rema}
 The no non-trivial expanding Ricci breather theorem was proved
 by Perelman in \cite{P} by discussing the
monotonicity of a scale invariant eigenvalue $\bar{\lambda}(\bg) =
\lambda(\bg) V^\frac{2}{n}(\bg)$. We present a different proof in
section \ref{Section No expanding breathers theorem}.
\end{rema}

\begin{rema}
For compact shrinking Ricci breathers, one can find references in
\cite{P}, \cite{BM}, \cite{BNL}, \cite{Hu}, \cite{I} and others.
\end{rema}

Throughout this paper, we use Einstein convention, i.e. repeated index implies summation.\\

\section{\textnormal{\bf Eigenvalues of $(-\Delta + \frac{R}{2})$ on manifolds}}\label{Section A remark on the assumption of positive Ricci curvature operator}

 Let ($M,\bg(t)$ be a compact Riemannian manifold, where $\bg(t)$ is a smooth solution to the Ricci flow equation on $0\le t<T$. In \cite{C}, Cao considers the problem from
 a viewpoint different from the entropy functional method.
 They study the eigenvalues $\lambda$ and eigenfunctions $f$ of the operator
 $-\Delta + \frac{R}{2}$ with the normalization $\Dint_m f^2dv = 1$. First, they define :
\begin{equation}\label{lowest defn1}
\lambda(h,t)=\Dint_M (-\Delta h + \frac{R}{2}h)h d\mu,
\end{equation}
 where $h$ is
a smooth function satisfying $\frac{d}{dt}(\Dint_M h^2 d\mu) = 0$
and $\Dint_M h^2 d\mu=1$.

They then obtain the monotonicity formula under the non-negative
curvature operator assumption as following:


\begin{theo}\label{steady positive thm1}\cite{C}
On a compact Riemannian manifold with nonnegative curvature
operator, the eigenvalues of the operator $-\Delta + \frac{R}{2}$
are nondecreasing under the Ricci flow, i.e.
\begin{equation} \label{steady positive equ1}
\frac{d}{dt}\lambda(f,t)|_{t=t_0} = 2\Dint R_{ij}f_if_j d\mu +
\Dint|Rc|^2f^2d\mu\ge 0.
\end{equation}
\end{theo}

In this theorem, $f_i$ denotes the covariant derivative of $f$ with
respect to $\frac{\partial}{\partial x_i}$, (also denoted as
$\partial_i$) and in (\ref{steady positive equ1}),
$\frac{d}{dt}\lambda(f,t)$ is evaluated at time $t=t_0$ and $f$ is
the corresponding eigenvalue at time $t_0$. As a direct consequence
of Theorem \ref{steady positive thm1}, they prove the following
\begin{coro}\label{steady
positive thm2}\cite{C} There is no compact steady Ricci breather
with nonnegative curvature operator, other than the one which is
Ricci-flat.
\end{coro}

\begin{rema}\label{eigenvalue rema1}
Clearly, at time $t$, if $f$ is the eigenfunction of the eigenvalue
$\lambda(t)$, then $\lambda(f,t) = \lambda(t)$.
\end{rema}

In this section, based on (\ref{steady positive equ1}), we will drop
the {\em curvature assumption} on the manifold and prove Theorem
\ref{thmljf1}.

\begin{proof}({\bf Theorem \ref{thmljf1}})
Let $\varphi$ be a function satisfying $f^2(x) = e^{-\varphi(x)}$
and plug it into (\ref{steady positive equ1}), we have
\begin{equation}\label{steady positive equ2}
\begin{array}{rll}
2\frac{d}{dt}\lambda(t)= 2\frac{d}{dt}\lambda(f,t)= \Dint
R_{ij}\nabla^i \varphi\nabla^j\varphi e^{-\varphi}d\mu + 2\Dint
|Rc|^2e^{-\varphi}d\mu.
\end{array}
\end{equation}

Using the divergence theorem and rearrangements, we derive the first
term of the last identity
\begin{equation}
\begin{array}{rll}
\Dint R_{ij}\nabla^i \varphi\nabla^j\varphi\ e^{-\varphi}d\mu
&=\Dint \nabla^i R_{ij}\nabla^j\varphi\ e^{-\varphi}d\mu+\Dint
R_{ij}\nabla^i \nabla^j\varphi\ e^{-\varphi}d\mu.\\
\end{array}
\label{steady positive equ3}
\end{equation}

By the contracted second Bianchi identity
``$\nabla^iR_{ij}=\frac{1}{2}\nabla_j R$'' and integration by parts,
we have
\begin{equation}
\begin{array}{rll}
(\ref{steady positive equ3}) &=\frac{1}{2}\Dint \nabla_j
R\nabla^j\varphi\ e^{-\varphi}d\mu+\Dint
R_{ij}\nabla^i \nabla^j\varphi\ e^{-\varphi}d\mu\\
&=\frac{1}{2}\Dint R\Delta e^{-\varphi}d\mu+\Dint R_{ij}\nabla^i
\nabla^j\varphi\ e^{-\varphi}d\mu.\\
\end{array}
\label{steady positive equ4}
\end{equation}
This also implies the following
\begin{equation}\label{steady positive equ4.1}
\Dint R_{ij}\nabla^i \nabla^j\varphi\ e^{-\varphi}d\mu = \Dint
R_{ij}\nabla^i \varphi\nabla^j\varphi\ e^{-\varphi}d\mu
-\frac{1}{2}\Dint R\Delta e^{-\varphi}d\mu.\\
\end{equation}
 On the other hand, using integration by parts and symmetry
of the hessian of functions, we have
\begin{equation}
\begin{array}{rll}
\Dint|\nabla\nabla\varphi|^2\ e^{-\varphi}d\mu &=
-\Dint\nabla_j\varphi \nabla_i \nabla^i \nabla^j\varphi\
e^{-\varphi}d\mu -\Dint\nabla_j\varphi \nabla^i \nabla^j\varphi
\nabla_i e^{-\varphi}d\mu\\
&= -\Dint\nabla_j\varphi\nabla_i \nabla^j \nabla^i\varphi\
e^{-\varphi}d\mu +\Dint\frac{1}{2}|\nabla\varphi|^2 \Delta
e^{-\varphi}d\mu.\\
\end{array}
 \label{steady positive equ5}
\end{equation}
By the commutator formulae for covariant derivatives which are known
as Ricci identities, see page 286 in \cite{BD}, we have
\begin{equation}
\nabla_i \nabla_j \nabla^i\varphi = \nabla_j \nabla_i
\nabla^i\varphi - R^k_{iji}\nabla^k \varphi = \nabla_j \nabla_i
\nabla^i\varphi + R_{kj}\nabla^k \varphi,
 \label{steady positive equ6}
\end{equation}
where $R^l_{ijk}$ represents the Riemann curvature $(3,1)$-tensor
and $R_{kj}$ denotes the Ricci curvature tensor. We use the
following convention of Riemann curvature tensor through out this
paper:
\[
R_m(X,Y)Z\equiv [\nabla_X,\nabla_Y]Z-\nabla_{[X,Y]}Z,
\]
where $\{X,Y,Z\}$ are vector fields on the manifold.

Combing (\ref{steady positive equ5}), (\ref{steady positive equ6}),
and the contracted second Bianchi identity, we have
\begin{equation} \label{steady positive equ7}
\begin{array}{rll}
\Dint|\nabla\nabla\varphi|^2& e^{-\varphi}d\mu =
-\Dint\nabla_j\varphi (\nabla_j \nabla^i \nabla^i\varphi -
R_{jk}\nabla^k \varphi)\ e^{-\varphi}d\mu+\Dint\frac{1}{2}|\nabla\varphi|^2 \Delta e^{-\varphi}d\mu\\
=& -\Dint\nabla_j\varphi \nabla^j \Delta\varphi\ e^{-\varphi}d\mu +
\Dint R_{jk}\nabla^j\varphi\nabla^k \varphi\ e^{-\varphi}d\mu
+\Dint\frac{1}{2}|\nabla\varphi|^2 \Delta e^{-\varphi}d\mu\\
=& -\Dint\Delta e^{-\varphi} \Delta\varphi d\mu + \Dint (\nabla_j
R_{jk}\nabla^k \varphi + R_{jk}\nabla^j\nabla^k \varphi)\
e^{-\varphi}d\mu\\
&\quad\quad  +\Dint\frac{1}{2}|\nabla\varphi|^2 \Delta
e^{-\varphi}d\mu\\
=& -\Dint\Delta e^{-\varphi} \Delta\varphi d\mu - \Dint
\frac{1}{2}\nabla_k R\nabla^k \varphi\ e^{-\varphi}d\mu -\Dint
R_{jk}\nabla^j\nabla^k \varphi\ e^{-\varphi}d\mu\\
&\quad\quad +\Dint\frac{1}{2}|\nabla\varphi|^2 \Delta e^{-\varphi}d\mu\\
=& -\Dint\Delta e^{-\varphi}\big( \Delta\varphi+\frac{1}{2}R -
\frac{1}{2}|\nabla\varphi|^2 \big)d\mu -\Dint R_{jk}\nabla^j\nabla^k
\varphi\ e^{-\varphi}d\mu.
\end{array}
\end{equation}
\indent We notice that we are free to change the dummy index from
$\{i,j,k,l\}$ to other index or exchange among them whenever
necessary.

Combing (\ref{steady positive equ4.1}) and (\ref{steady positive
equ7}), we have the following
\begin{equation}
\begin{array}{rll}
\Dint 2R_{ij}\nabla^i\nabla^j \varphi e^{-\varphi}d\mu &+ \Dint
|\nabla\nabla \varphi|^2 e^{-\varphi}d\mu \\
= &\Dint R_{jk}\nabla^j\nabla^k \varphi\
e^{-\varphi}d\mu-\Dint\Delta e^{-\varphi}\big(
\Delta\varphi+\frac{1}{2}R
-\frac{1}{2}|\nabla\varphi|^2 \big)d\mu\\
= &\Dint R_{ij}\nabla^i \varphi\nabla^j\varphi\ e^{-\varphi}d\mu
-\Dint\Delta e^{-\varphi}\big( \Delta\varphi+R
-\frac{1}{2}|\nabla\varphi|^2 \big)d\mu\\
\end{array}
 \label{steady positive equ8}
\end{equation}

\indent Recall (\ref{thmljf1 rev equ1}), then a simple calculation
yields
\begin{equation}\label{steady positive equ9}
2\lambda(t) = \Delta_{\bg(t)}\varphi+R_{\bg(t)}
-\frac{1}{2}|\nabla\varphi|_{\bg(t)}^2
\end{equation}

Plugging (\ref{steady positive equ9}) into (\ref{steady positive
equ8}), by divergence theorem on closed manifolds, we have
\begin{equation}\label{steady positive equ10}
\begin{array}{rll}
\Dint 2R_{ij}\nabla^i\nabla^j \varphi e^{-\varphi}d\mu &+ \Dint
|\nabla\nabla \varphi|^2 e^{-\varphi}d\mu \\
=&\Dint R_{ij}\nabla^i\varphi \nabla^j\varphi\ e^{-\varphi}d\mu
-\Dint 2\lambda(t)\Delta e^{-\varphi}d\mu\\
=&\Dint R_{ij}\nabla^i \varphi\nabla^j\varphi\ e^{-\varphi}d\mu.
\end{array}
\end{equation}
In the end, we plug (\ref{steady positive equ10}) into (\ref{steady
positive equ2}) and have the following
\begin{equation}\label{steady positive equ11}
\begin{array}{rll}
2\frac{d}{dt}\lambda(t) &= \Dint R_{ij}\nabla^i
\varphi\nabla^j\varphi e^{-\varphi}d\mu + \Dint
|Rc|^2e^{-\varphi}d\mu + \Dint
|Rc|^2e^{-\varphi}d\mu\\
&= \Dint |R_{ij}+\nabla^i\nabla^j \varphi|^2\ e^{-\varphi}d\mu +
\Dint
|Rc|^2\ e^{-\varphi}d\mu\\
&\ge 0
\end{array}
\end{equation}
\end{proof}


\begin{rema}
One can find applications of Theorem \ref{thmljf1} in \cite{C} {\rm
(} e.g. Theorem $3${\rm)}. When the lowest eigenvalues are
concerned, the eigenfunctions are always smooth and positive, see in
\cite{BM}. For the monotonicity of eigenvalues of ordinary Laplace
operator under Ricci flow, one can see a result in \cite{M}.
\end{rema}

\begin{rema}
For evolution of Yamabe constant under Ricci flow, see a recent
preprint \cite{CL}. We thank Xiao-Dong Cao for pointing out the
reference to us.
\end{rema}
%
%


\section{\textnormal{\bf Construct entropy functional from first variation formula}}\label{Section From first variation formula to Entropy functionals}
In this section, we will construct a functional $\mathcal{E}$ and
show that it is monotone along the Ricci flow.

Under the following coupled system:
\begin{equation}\label{from equ1}
\left\{
\begin{array}{rll}
&\frac{\partial}{\partial t}\mathbf{g}_{ij} = -2(R_{ij}+\nabla_i\nabla_j f)\\
&\frac{\partial}{\partial t}f = -\Delta f - R,
\end{array}
\right.
\end{equation}
one has the following first variation formula, (see \cite{P})
\begin{equation}\label{from equ2}
\frac{d}{d t}\mathcal{F} = 2\Dint_M |R_{ij} + \nabla_i\nabla_j f|^2
e^{-f}d\mu \ge 0,
\end{equation}
where $\cF$ is defined as in (\ref{Intro def1}).  Hence, the
modified Ricci flow can be viewed as an $L^2$ gradient Ricci flow of
Perelman's $\mathcal{F}$-functional.

One natural question to ask is : can we find a functional
$\mathcal{E}$ such that the `honest' Ricci flow is the $L^2$
gradient Ricci flow of it in a certain sense? Namely, one expects
that under the following coupled system :
\begin{equation}\label{from equ3}
\left\{
\begin{array}{rll}
&\frac{\partial}{\partial t}\mathbf{g}_{ij} = -2R_{ij}\\
&\frac{\partial}{\partial t}f = -\Delta f - R + |\nabla f|^2,
\end{array}
\right.
\end{equation}
there is a functional $\cE$ such that the following is true,
\begin{equation}\label{from equ4}
\frac{d}{d t}\mathcal{E} = 2\Dint_M |Rc|^2 e^{-f}d\mu \ge 0.
\end{equation}
It turns out that this can be done by integrating the above first
variation formula over a time interval $[0,t]$. We carry out the
computations as the following
\begin{equation}\label{from equ5}
\begin{array}{rll}
\mathcal{E}(t) - \mathcal{E}(0) &= \Dint_0^t \frac{d}{d
s}\mathcal{E} ds = 2\Dint_0^t \Dint_M
|Rc|^2 e^{-f}d\mu ds\\
& =2\Dint_M \Dint_0^t \mathbf{g}^{ij}\mathbf{g}^{kl}R_{ik}R_{jl}\ e^{-f}\sqrt{{\rm det\ }\mathbf{g}}\ ds\ dx^n\\
& =\Dint_M \Dint_0^t \frac{\partial\mathbf{g}^{jl}}{\partial t}R_{jl}\ e^{-f}\sqrt{{\rm det\ }\mathbf{g}}\ ds\ dx^n\\
\end{array}
\end{equation}
The last step follows from a lemma in \cite{BD}, see page 67, {Lemma
3.1}.

Integration by parts yields
\begin{equation}\label{from equ6}
\begin{array}{rll}
\mathcal{E}(t) - \mathcal{E}(0) & =\Dint_M
\bigg[\mathbf{g}^{jl}R_{jl} e^{-f}\sqrt{{\rm det\
}\mathbf{g}}\bigg|_0^s-\Dint_0^t \mathbf{g}^{jl}\frac{\partial
R_{jl}}{\partial t} e^{-f}\sqrt{{\rm det\ }\mathbf{g}}\ ds\\
&\quad-\Dint_0^t \mathbf{g}^{jl}R_{jl}\frac{\partial
e^{-f}\sqrt{{\rm det\ }\mathbf{g}}}{\partial t} \ ds\bigg]\ dx^n\\
& =\Dint_M R e^{-f}d\mu\bigg|_0^s-\Dint_M \bigg[\Dint_0^t
\mathbf{g}^{jl}\frac{\partial
R_{jl}}{\partial t} e^{-f}\sqrt{{\rm det\ }\mathbf{g}}\ ds\\
&\quad+\Dint_0^t \mathbf{g}^{jl}R_{jl}\frac{\partial
e^{-f}\sqrt{{\rm det\ }\mathbf{g}}}{\partial t} \ ds\bigg]\
dx^n\\\end{array}
\end{equation}
Using {Lemma 3.5} in \cite{BD} and the fact that metric tensor
$\mathbf{g}$ is a covariant constant, by letting $h=-2Rc$, we have
the following \footnote{ As suggested by the referee, a simpler way
to obtain equation (\ref{from equ7}) is :
\[\begin{array}{rll}\frac{\partial R}{\partial t} =& \Delta R + 2 |Rc|^2=
\mathbf{g}^{ij}\frac{\partial R_{ij}}{\partial t} +
R_{ij}\frac{\partial \mathbf{g}^{ij}}{\partial t}\\
=&\mathbf{g}^{ij}\frac{\partial R_{ij}}{\partial t} + 2
|Rc|^2.\end{array}
\] }
\begin{equation}\label{from equ7}
\begin{array}{rll}
\displaystyle\mathbf{g}^{jl}\frac{\partial R_{jl}}{\partial t} &=
-\mathbf{g}^{jl}\mathbf{g}^{pq}(\nabla_q\nabla_j
R_{lp}+\nabla_q\nabla_l
R_{jp}-\nabla_q\nabla_p R_{jl}-\nabla_j\nabla_l R_{pq})\\
&=\Delta R.
\end{array}
\end{equation}
\indent In the above computations, we have used the contracted
second Bianchi identity. With the help of {Lemma 3.9} in \cite{BD},
We obtain
\begin{equation}\label{from equ8}
\begin{array}{rll}
\displaystyle\frac{\partial e^{-f}\sqrt{{\rm det\
}\mathbf{g}}}{\partial t}&= (-\frac{\partial}{\partial
t}f-R)e^{-f}\sqrt{{\rm det\ }\mathbf{g}}
\end{array}
\end{equation}
Plug (\ref{from equ7}) and (\ref{from equ8}) into (\ref{from equ6}),
we have
\begin{equation}\label{from equ9}
\begin{array}{rll}
\mathcal{E}(t) - \mathcal{E}(0) &  =\Dint_M R
e^{-f}d\mu\bigg|_0^s-\Dint_0^t\Dint_M \bigg[ \Delta R\ e^{-f} +R
(-\frac{\partial}{\partial
t}f-R)e^{-f}\bigg]\ d\mu\ ds\\
&  =\Dint_M R e^{-f}d\mu\bigg|_0^s-\Dint_0^t\Dint_M R \bigg[ -\Delta
f +|\nabla f|^2
 -\frac{\partial}{\partial
t}f-R\bigg ] e^{-f}\ d\mu\ ds\\
&  =\Dint_M R e^{-f}d\mu\bigg|_{s\ =\ t}-\Dint_M R
e^{-f}d\mu\bigg|_{s\ =\ 0}.
\end{array}
\end{equation}
The last step follows from the third equation in the coupled system
(\ref{from equ3}). The computation above suggests that we could
define $\mathcal{E}$-functional as
\begin{equation}\label{from equ10}
\mathcal{E} = \Dint_M R e^{-f}d\mu.
\end{equation}

Next, we prove the monotonicity property for
$\mathcal{E}$-functional directly.

\begin{prop}
Assume $\mathbf{g}(t)$ satisfies the Ricci flow equation over the
time interval [0,T), and also function $f$ satisfies the evolution
equation (\ref{sec1.1}) (also (\ref{from equ3})), then
\[
\frac{d\mathcal{E}}{d t} = 2\Dint_M |Rc|^2e^{-f}d\mu.
\]
\label{prop1}
\end{prop}

\begin{proof}({\bf First proof of Proposition \ref{prop1}})
 Using $\frac{\partial R}{\partial t} = \Delta R + 2|Rc|^2$, we derive
the following
\[
\begin{array}{rll}
\frac{d}{d t}\Dint_M Re^{-f}d\mu =& \Dint_M \frac{dR}{d t}e^{-f}d\mu
+
R\frac{d}{d t}(e^{-f}d\mu)\\
=& \Dint_M (\Delta R +
2|Rc|^2)e^{-f}d\mu +\Dint_M R(-f_t - R)e^{-f}d\mu\\
=& 2\Dint_M|Rc|^2e^{-f}d\mu + \Dint_M R [-\Delta f + |\nabla f|^2-f_t - R] e^{-f} d\mu \\
=& 2\Dint_M|Rc|^2e^{-f}d\mu.
\end{array}
\label{2.3}
\]

The last equality comes from the second equation of the coupled
system (\ref{sec1.1}).

({\bf Second proof}) Under the following coupled system with
modified Ricci flow (see also similar system in (\ref{Entropy
formulae over Expanders thm1 equ1})),
\begin{equation}\label{prop1 2}
\left\{
\begin{array}{rll}
\Dfrac{\partial}{\partial t}\mathbf{g}_{ij} &=
-2\bigg(R_{ij}+\nabla_i\nabla_j f\bigg)\\
\frac{\partial}{\partial t}f &= -\Delta f - R,
\end{array} \right.
\end{equation}
we have the following
\[
\begin{array}{rll}
\frac{d}{d t}\Dint_M Re^{-f}d\mu = &\Dint_M (\Delta R -\nabla
R\nabla f + 2|Rc|^2) e^{-f}\ d\mu\\
=&2\Dint_M |Rc|^2 e^{-f}\ d\mu.
\end{array}
\label{2.4}
\]
By using the diffeomorphism invariance, see e.g. Proposition 1.2 in
\cite{P}, we prove that under the original coupled system we still
have
\[
\begin{array}{rll}
\frac{d}{d t}\Dint_M Re^{-f}d\mu = 2\Dint_M |Rc|^2 e^{-f}\ d\mu.
\end{array}
\]

\end{proof}

\begin{rema}
We will apply the diffeomorphism invariance principle generally in
this paper. See details in the proof of Corollary \ref{Entropy
formulae over Expanders cor1} and Corollary \ref{Entropy formulae
over Expanders cor2}.
\end{rema}

\section{\textnormal{\bf Entropy functionals $\mathcal{F}_k$ and their monotonicity}}\label{Section Entropy functionals}

\begin{defi}
We define the following variations of $\mathcal{F}$-functional,
\begin{equation}
\mathcal{F}_k(\mathbf{g},f) = \Dint_M (k R + |\nabla
f|^2)e^{-f}d\mu,\label{2.2}
\end{equation}
where $k\ge 1$. When $k=1$, this is the $\mathcal{F}$-functional.
\end{defi}

Next we derive the monotonicity formula for these functionals
$\mathcal{F}_k(\mathbf{g},f)$.

\begin{theo}\label{thm2}
Suppose the Ricci flow of $\mathbf{g}(t)$ exists for $[0,T)$, then
all the functionals $\mathcal{F}_k(\mathbf{g},f)$ will be monotone
under the coupled system (\ref{from equ3}), i.e.
\begin{equation}\label{thm2 equ1}
\begin{array}{rll}
\frac{d}{d t}\mathcal{F}_k(\mathbf{g}_{ij},f)=2(k-1)\Dint_M
|Rc|^2e^{-f}d\mu+2\Dint_M |R_{ij}+\nabla_i\nabla_j f|^2e^{-f}d\mu
\ge 0
\end{array}
\end{equation}
Furthermore, the monotonicity is strict unless the Ricci flow is a
trivial Ricci soliton and $f$ is a constant function. Namely, the
metric is Ricci-flat.
\end{theo}
\begin{proof}({\bf Theorem \ref{thm2}})
Under the coupled system (\ref{sec1.1})
\[
\begin{cases}
 &\frac{\partial}{\partial t}\mathbf{g}_{ij} = -2R_{ij}\\
&\frac{\partial}{\partial t}f = -\Delta f - R + |\nabla f|^2,
\end{cases}
\]
we have shown in {Proposition \ref{prop1}} that,
\[
\frac{d}{d t}\Dint_M Re^{-f}d\mu = 2\Dint_M |Rc|^2e^{-f}d\mu.
\]
On the other hand, in \cite{P}, it was shown that under the same
system (\ref{sec1.1}),
\[
\frac{d}{d t}\Dint_M (R+|\nabla f|^2)e^{-f}d\mu = 2\Dint_M
|R_{ij}+\nabla_i\nabla_j f|^2e^{-f}d\mu.
\]
Put them together we get the following formula:
\begin{equation}
\begin{array}{rll}
\frac{d}{d t}\mathcal{F}_k(\mathbf{g},f)&=\frac{d}{d t}\Dint_M (kR+|\nabla f|^2)e^{-f}d\mu\\
 &=2(k-1)\Dint_M
|Rc|^2e^{-f}d\mu+2\Dint_M |R_{ij}+\nabla_i\nabla_j
f|^2e^{-f}d\mu\\
&\ge 0.
\end{array}\label{equ2.4}
\end{equation}
For $k=1$, the above is consistent with the monotonicity formula of
$\mathcal{F}$-functional. The second statement is obvious. This
finishes the proof of {\em Theorem \ref{thm2}}.
\end{proof}

\begin{rema}
We notice that under the coupled system (\ref{sec1.1}), the Ricci
flow can be viewed as a $\mathbb{L}^2$ gradient flow of Perelman's
$\mathcal{F}$ functional up to a diffeomorphism where our
functionals are not. But the monotonicity is still retained.
\end{rema}

\begin{rema}
 Our functionals yield more information about the Ricci
tensor itself.
\end{rema}

\section{\textnormal{\bf Eigenvalues and compact steady Ricci breathers}}\label{section steady}
 In this section, we discuss the applications
of the monotonicity formula (\ref{thm2 equ1}) we derived in Theorem
\ref{thm2}.

First, we recall the definition of Ricci breathers, see original
definition in \cite{I} and \cite{P}.
\begin{defi}\label{breather compact def1}
A metric $\mathbf{g}(t)$ evolving by the Ricci flow is called a
breather, if for some $t_1<t_2$ and $\alpha >0$ the metrics $\alpha
\mathbf{g}(t_1)$ and $\mathbf{g}(t_2)$ differ only by a
diffeomorphism; the cases $\alpha =1 $, $\alpha <1 $ , $\alpha >1 $
correspond to steady, shrinking or expanding breathers,
respectively.
\end{defi}
Trivial breathers are called Ricci solitons for which the above
properties are true for each pair of $t_1$ and $t_2$.

Define $\lambda_k(\mathbf{g}) = \inf \mathcal{F}_k(\mathbf{g},f)$,
where infimum is taken over all smooth $f$, satisfying $\int_M
e^{-f}d\mu = 1$. $\lambda_k$ is the lowest eigenvalue of the
corresponding operators $-4\Delta + kR$ for $k>1$. By applying
direct methods and elliptic regularity theory (see
\cite{GT},\S8.12), one can see that the infimum is always attained.

Using the monotonicity in Theorem \ref{thm2}, we have

\begin{theo}\label{revise thmljf1}
On a compact Riemannian manifold ($M,\bg(t)$), where $\bg(t)$
satisfies the Ricci flow equation for $t\in [0,T)$, the lowest
eigenvalue $\lambda_k$ of the operator $-4\Delta + kR$ is
nondecreasing under the Ricci flow. The monotonicity is strict
unless the metric is Ricci-flat.
\end{theo}

\begin{proof}({\bf Theorem \ref{revise thmljf1}})
For any time $t_1$, $t_2$ in $[0,T)$, suppose that at time $t_2$,
the lowest eigenvalue $\lambda_k(\mathbf{g}(t_2))$ is attained by a
function $f_k(x)$. Evolving under the backward Ricci flow, we get a
solution $f_k(x,t)$ to the coupled system (\ref{sec1.1}) which
satisfies the initial condition $f_k(x,t_2)=f_k(x)$.

Using the monotonicity formula of (\ref{thm2 equ1}), we have
$$
\begin{array}{rll}
\lambda_k(\mathbf{g}(t_2)) &= \mathcal{F}_k(\mathbf{g}(t_2),f_k(t_2))\\
&\ge \mathcal{F}_k(\mathbf{g}(t_1),f_k(t_1))\\
&\ge \inf \mathcal{F}_k(\mathbf{g}(t_1),f)\\
&= \lambda_k(\mathbf{g}(t_1)).
\end{array}
$$

This proves that $\lambda_k$ is nondecreasing under the Ricci flow.
Since the monotonicity of $\cF_k$ is strict unless the metric is
Ricci-flat, this finishes the proof of Theorem \ref{revise thmljf1}.
\end{proof}

As a corollary, we generalize the theorem of Cao in the case of the
lowest eigenvalue.

\begin{coro}\label{Cao thm1}
On a compact Riemannian manifold, the lowest eigenvalues of the
operator $-\Delta + \frac{R}{2}$ are nondecreasing under the Ricci
flow.
\end{coro}

\begin{proof}({\bf Corollary \ref{Cao thm1}})
Let $k=2$, then $\frac{1}{4}\lambda_2$ is the lowest eigenvalue of
$-\Delta + \frac{R}{2}$ and the result will follow.
\end{proof}


As an application of Theorem \ref{steady positive thm1}, Cao proved
the following

\begin{theo}\label{steady
positive thm2}\cite{C} There is no compact steady Ricci breather
with {\em nonnegative curvature operator}, other than the one which
is Ricci-flat.
\end{theo}

As a corollary of Theorem \ref{revise thmljf1}, we drop the {\em
nonnegative curvature operator} condition  and have the following
\begin{coro}\label{thmljf2}
There is no compact steady Ricci breather other than the one which
is Ricci-flat.
\end{coro}

\begin{proof}({\bf Theorem \ref{thmljf2}})
For a Ricci breather, let $t_1$ and $t_2$ be as in the definition,
then $\lambda_k(t_1)=\lambda_k(t_2)$ for a steady breather due to
the diffeomorphism invariance of the eigenvalues. The fact that
$\lambda_k(t)$ fails to be strictly increasing yields that the
manifold must be Ricci-flat.
\end{proof}

\begin{rema}\label{chap1 intro rmk2}
This result was first proved by T. Ivey in \cite{I} with a different
approach. See also details in the book \cite{BM}.
\end{rema}

\section{\textnormal{\bf New formulae over expanders }}\label{Section New entropy formulae over expanders}

In this section, we define the following functionals and discuss
their first variation formulae under modified Ricci flow and Ricci
flow respectively
\begin{equation}\label{Entropy formulae over Expanders equ1}
\begin{array}{rll}
\mathcal{W}_e(\mathbf{g},\tau(t),f)&=\tau^2\Dint_M \big[R+\frac{n}{2\tau}+\Delta f\big]e^{-f}d\mu,\\
\cW_{ek}(\mathbf{g},\tau(t),f)&=\tau^2\Dint_M
\big[k(R+\frac{n}{2\tau})+\Delta f\big]e^{-f}d\mu,
\end{array}
\end{equation}
where $k > 1$.


\begin{rema} \label{New entropy formulae over expanders rmk1}
The functionals we obtain in this section are different from
Perelman's $\cW$-functional. They are not scale invariant.
\end{rema}

\begin{rema} \label{New entropy formulae over expanders rmk2}
After this paper was submitted, we found out that in \cite{P},
Perelman has defined a functional similar to our functional
$\mathcal{W}_e$ which he called $-<E>$. Perelman's $<E>$ is modeled
on shrinking Ricci solitons and ours are the corresponding version
on expanding solitons. Furthermore, the functionals
$\mathcal{W}_{ek}$ are new to our knowledge. A different motivation
which made us discover these functionals, $\mathcal{W}_e$ and
$\mathcal{W}_{ek}$, will appear somewhere else.
\end{rema}

\begin{rema} \label{New entropy formulae over expanders rmk3}
In \cite{FIN}, M. Feldman, T. Ilmanen, and L. Ni constructed a scale
invariant $\mathcal{W}$-entropy which is an analogue of Perelman's
$\mathcal{W}$-entropy but has vanishing first variation over
expanders. There is also a very good unified treatment about entropy
formulae over steady, expanding, and shrinking Ricci breathers in
\cite{BM}. See a related generalization of Perelman¡¯s formula also
in a recent preprint \cite{KZ}.
\end{rema}


We start with deriving the first variation formulae of
$\mathcal{W}_e$ and $\cW_{ek}$.

\begin{theo}\label{Entropy formulae over Expanders thm1}
Under the following coupled system
\begin{equation}\label{Entropy formulae over Expanders thm1 equ1}
\left\{
\begin{array}{rll}
\Dfrac{\partial}{\partial t}\mathbf{g}_{ij} &=
-2\bigg(R_{ij}+\nabla_i\nabla_j f\bigg)\\
\frac{\partial}{\partial t}f &= -\Delta f - R\\
\Dfrac{d}{dt}\tau &= 1,
\end{array}
\right.
\end{equation}
the first variation formula for
$\mathcal{W}_e(\mathbf{g},\tau(t),f)$ is
\begin{equation}\label{Entropy formulae over Expanders thm1 equ2}
\begin{array}{rll}
\Dfrac{d}{d t} \mathcal{W}_e = 2\tau^2\Dint_M
\big|R_{ij}+\nabla_i\nabla_j
f+\Dfrac{1}{2\tau}\mathbf{g}_{ij}\big|^2e^{-f}d\mu.
\end{array}
\end{equation}
\end{theo}

\begin{theo}\label{Entropy formulae over Expanders thm2}
Under the coupled system (\ref{Entropy formulae over Expanders thm1
equ1}), the first variation formula for
$\cW_{ek}(\mathbf{g},\tau(t),f)$ is
\begin{equation}\label{Entropy formulae over Expanders thm2 equ2}
\begin{array}{rll}
\Dfrac{d}{d t} \cW_{ek} =& 2\tau^2\Dint_M
\big|R_{ij}+\nabla_i\nabla_j
f+\Dfrac{1}{2\tau}\mathbf{g}_{ij}\big|^2e^{-f}d\mu \\
&+ 2(k-1)\tau^2\Dint_M
\big|R_{ij}+\Dfrac{1}{2\tau}\mathbf{g}_{ij}\big|^2e^{-f}d\mu.
\end{array}
\end{equation}
where $k>1$.
\end{theo}

\begin{proof}({\bf Theorem \ref{Entropy formulae over Expanders thm1}})
The proof is by direct computations. We notice that $\mathcal{W}_e =
\tau^2 \mathcal{F} + \frac{n}{2\tau}$. Under the coupled system
(\ref{Entropy formulae over Expanders thm1 equ1}), it is known that
$\frac{d}{dt}\mathcal{F} = 2\Dint_M|R_{ij}+\nabla_i\nabla_j
f|^2e^{-f}d\mu$. This implies
\begin{equation}\label{Entropy formulae over Expanders thm1 proof equ1}
\begin{array}{rll}
\frac{d}{dt}\mathcal{W}_e &= 2\tau^2\Dint_M|R_{ij}+\nabla_i\nabla_j
f|^2e^{-f}d\mu
+ 2\tau\mathcal{F} + \frac{n}{2}\\
&=2\tau^2\Dint_M|R_{ij}+\nabla_i\nabla_j f +
\frac{1}{2\tau}\mathbf{g}_{ij}|^2e^{-f}d\mu.
\end{array}
\end{equation}
\end{proof}

\begin{proof}({\bf Theorem \ref{Entropy formulae over Expanders thm2}})
By definition, we have $\cW_{ek} = \mathcal{W}_e +
(k-1)\tau^2\Dint_M (R + \frac{n}{2\tau})e^{-f}d\mu$. Recall we
defined $\mathcal{E} = \Dint_M Re^{-f}d\mu$ in (\ref{from equ10}).
This yields $\cW_{ek} = \mathcal{W}_e + (k-1)(\tau^2\mathcal{E} +
\frac{n}{2}\tau)$.

As in the second proof of Proposition \ref{prop1}, we have, under
the coupled system (\ref{Entropy formulae over Expanders thm1
equ1}), $\frac{d}{dt}\mathcal{E} = 2\Dint_M |Rc|^2e^{-f}d\mu$. Using
Theorem \ref{Entropy formulae over Expanders thm1}, the rest of the
proof is direct computations as the following
\begin{equation}\label{Entropy formulae over Expanders thm2 proof equ1}
\begin{array}{rll}
\frac{d}{dt}\cW_{ek} =& 2\tau^2\Dint_M|R_{ij}+\nabla_i\nabla_j f +
\frac{1}{2\tau}\mathbf{g}_{ij}|^2e^{-f}d\mu \\
&\quad\quad+ (k-1)\bigg[2\tau^2\Dint_M|Rc|^2e^{-f}d\mu + 2\tau\mathcal{E} + \frac{n}{2}\bigg]\\
=& 2\tau^2\Dint_M|R_{ij}+\nabla_i\nabla_j f +
\frac{1}{2\tau}\mathbf{g}_{ij}|^2e^{-f}d\mu \\
&\quad\quad+ 2(k-1)\tau^2\int_M|R_{ij}+
\frac{1}{2\tau}\mathbf{g}_{ij}|^2e^{-f}d\mu.
\end{array}
\end{equation}
\end{proof}

As direct corollaries of Theorem \ref{Entropy formulae over
Expanders thm1} and Theorem \ref{Entropy formulae over Expanders
thm2}, we also obtain first variation formulae under the `honest'
Ricci flow.

\begin{coro}\label{Entropy formulae over Expanders cor1}
Under the following coupled system
\begin{equation}\label{Entropy formulae over Expanders cor1 equ1}
\left\{
\begin{array}{rll}
\Dfrac{\partial}{\partial t}\mathbf{g}_{ij} &= -2R_{ij}\\
\frac{\partial}{\partial t}f &= -\Delta f +|\nabla f|^2- R\\
\Dfrac{d}{dt}\tau &= 1,
\end{array}
\right.
\end{equation}
the first variation formula for
$\mathcal{W}_e(\mathbf{g},\tau(t),f)$
\begin{equation}\label{Entropy formulae over Expanders cor1 equ2}
\begin{array}{rll}
\Dfrac{d}{d t} \mathcal{W}_e = 2\tau^2\Dint_M
\big|R_{ij}+\nabla_i\nabla_j
f+\Dfrac{1}{2\tau}\mathbf{g}_{ij}\big|^2e^{-f}d\mu.
\end{array}
\end{equation}

\end{coro}

\begin{coro}\label{Entropy formulae over Expanders cor2}
Under the coupled system (\ref{Entropy formulae over Expanders cor1
equ1}), the first variation formula for
$\cW_{ek}(\mathbf{g},\tau(t),f)$ is
\begin{equation}\label{Entropy formulae over Expanders cor2 equ2}
\begin{array}{rll}
\Dfrac{d}{d t} \cW_{ek} =& 2\tau^2\Dint_M
\big|R_{ij}+\nabla_i\nabla_j
f+\Dfrac{1}{2\tau}\mathbf{g}_{ij}\big|^2e^{-f}d\mu \\
&\quad+ 2(k-1)\tau^2\Dint_M
\big|R_{ij}+\Dfrac{1}{2\tau}\mathbf{g}_{ij}\big|^2e^{-f}d\mu,
\end{array}
\end{equation}
where $k>1$.
\end{coro}

\begin{proof}({\bf Corollary \ref{Entropy formulae over Expanders cor1}})
We adapt a proof similar to the one of Lemma 5.15 in \cite{BM} (see
also in Proposition 1.2 of \cite{P}). We observe that the modified
Ricci flow and Ricci flow differ only by a diffeomorphism. Since
$(\mathbf{g}(t), f(t))$ is a solution to (\ref{Entropy formulae over
Expanders cor1 equ1}), the pair $(\bar{\mathbf{g}}(t), \bar{f}(t))$
defined by $\bar{\mathbf{g}}(t) = \Phi^*(t)\mathbf{g}(t)$ and
$\bar{f}(t) = f(t)\circ \Phi(t)$, is a solution to (\ref{Entropy
formulae over Expanders thm1 equ1}), where $\Phi(t):M \rightarrow M$
is a one parameter family of diffeomorphisms defined by

\begeq
\begin{array}{rll}
\frac{d}{dt}\Phi(t) &= \nabla_{\bar{\bg}(t)}\bar{f}(t) \\
\Phi(0)&= id_M.
\end{array}
\endeq
 Now
$\mathcal{W}_e(\mathbf{g},f)=\mathcal{W}_e(\bar{\mathbf{g}},\bar{f})$,
so that by (\ref{Entropy formulae over Expanders thm1 equ2}) we have
\begin{equation}\label{Entropy formulae over Expanders cor1 proof equ1}
\begin{array}{rll}
\Dfrac{d}{dt} \mathcal{W}_e(\mathbf{g}(t), f(t))
&=\Dfrac{d}{dt} \mathcal{W}_e(\bar{\mathbf{g}}(t), \bar{f}(t))\\
& = 2\tau^2\Dint_M \big|\bar{R}_{ij}+\bar{\nabla}_i\bar{\nabla}_j
f-\Dfrac{1}{2\tau}\bar{\mathbf{g}}_{ij}\big|^2e^{-\bar{f}}d\bar{\mu}\\
& = 2\tau^2\Dint_M \big|R_{ij}+\nabla_i\nabla_j
f-\Dfrac{1}{2\tau}\mathbf{g}_{ij}\big|^2e^{-f}d\mu.
\end{array}
\end{equation}
\end{proof}

\begin{proof}({\bf Corollary \ref{Entropy formulae over Expanders cor2}})
The proof is very similar to the proof of Corollary \ref{Entropy
formulae over Expanders cor1} and shall be skipped here.
\end{proof}

If we unify Corollary \ref{Entropy formulae over Expanders cor1} and
Corollary \ref{Entropy formulae over Expanders cor2}, we obtain
Theorem \ref{Entropy formulae over Expanders main} given in the
introduction section.

\begin{rema}
One can define exactly similar new formulae $\mathcal{W}_s$ and
$\mathcal{W}_{sk}$ over shrinkers which differ only by replacing
$\frac{n}{2\tau}$ by $-\frac{n}{2\tau}$ and letting
$\frac{d}{dt}\tau = -1$. Similar monotonicity formulae still hold.
\end{rema}


\section{\textnormal{\bf No expanding breathers theorem}}\label{Section No expanding breathers theorem}

We will focus on compact expanding Ricci breathers. First, we define
$\mu_e(\mathbf{g},\tau) = \inf {\mathcal W}_e(\mathbf{g},\tau,f)
=\inf \tau^2\Dint_M \big[R+\frac{n}{2\tau}+\Delta f\big]e^{-f}d\mu
$, where the infimum is taken over all the smooth functions
satisfying $\Dint_Me^{-f}d\mu = 1$.


We prove the following result of no non-trivial expanding Ricci
breathers which was obtained by Perelman in \cite{P} with a
different method.

\begin{coro}\label{No expanding breathers theorem thm1}
There is no expanding Ricci breather on compact Riemannian manifolds
other than expanding gradient Ricci solitons.
\end{coro}

\begin{proof}({\bf Corollary \ref{No expanding breathers theorem thm1}})
By the definition of Ricci breathers, there exist a pair of time
moments $t_1$ and $t_2$, such that the Ricci flow solution
$\mathbf{g}(t)$ at these two moments differ only by a diffeomorphism
and a scaling $\alpha$, i.e., $\mathbf{g}(t_2) = \alpha\Phi^*
\mathbf{g}(t_1)$, where $\Phi$ is a diffeomorphism and the scalar
$\alpha >1$.

By the standard argument about the existence of a minimizer of
$\mu_e$ at a fixed time moment, there exists a smooth function $f$
which attains the infimum. By solving the backward heat equation in
Corollary \ref{Entropy formulae over Expanders cor1}, we get a
smooth function $f(x,t): M\times [t_1, t_2]\longrightarrow
\mathbb{R}$ with initial condition $f(\cdot,t_2) = f(\cdot)$. We
define a linear function $\tau:[t_1,t_2]\longrightarrow (0,+\infty)$
with $\tau(t_2)=T+t_2$, where $T$ is a real number. Under the
coupled system (\ref{Entropy formulae over Expanders cor1 equ1}), by
the monotonicity formula, we obtain the following
\begin{equation}\label{No expanding breathers theorem thm1 proof equ2}
\begin{array}{rll}
\mu_e(\mathbf{g}(t_1),\tau(t_1)) &= \inf\big|_{t = t_1}\mathcal{W}_e(\mathbf{g},\tau,f)\\
& \le \mathcal{W}_e(\mathbf{g},\tau,f)\big|_{t = t_1}\\
& \le \mathcal{W}_e(\mathbf{g},\tau,f)\big|_{t = t_2}\\
&=\mu_e(\mathbf{g}(t_2),\tau(t_2)).
\end{array}
\end{equation}
On the other hand, it is easy to see that, if we simultaneously
scale $\tau$ and $\mathbf{g}$ by a scalar $\alpha$, we have the
following scaling property of $\mu_e$,
\begin{equation}\label{No expanding breathers theorem thm1 proof equ3}
\mu_e(\alpha\mathbf{g}(t),\alpha\tau) =
\alpha\mu_e(\mathbf{g}(t),\tau).
\end{equation}
Combining (\ref{No expanding breathers theorem thm1 proof equ2}),
(\ref{No expanding breathers theorem thm1 proof equ3}), and the
diffeomorphic invariant property of the functionals, we have
\begin{equation}\label{No expanding breathers theorem thm1 proof equ4}
\mu_e(\mathbf{g}(t_1),\tau(t_1)\le \mu_e(\mathbf{g}(t_2),\tau(t_2))
= \alpha\mu_e(\mathbf{g}(t_1),\tau(t_1)).
\end{equation}
This yields $0\le (\alpha -1)\mu_e(\mathbf{g}(t_1),\tau(t_1))$, i.e.
\begin{equation}\label{No expanding breathers theorem thm1 proof equ5}
\mu_e(\mathbf{g}(t_1),\tau(t_1))\ge 0.
\end{equation}

Since $\mathbf{g}(t_2)$ and $\mathbf{g}(t_1)$ are differed by the
scalar $\alpha$ up to a diffeomorphism, in order to scale
$\mathbf{g}$ and $\tau$ simultaneously at time $t_2$, we need to
make $\tau(t_2) = \alpha \tau(t_1)$. This yields that $T+t_2 =
\alpha(T+t_1)$. Equivalently, we can choose $T = \frac{t_2-\alpha
t_1}{\alpha -1}$. \footnote{ By the choice of $T$, we can show that
$ \tau(t_1) =T+t_1 =\frac{t_2-t_1}{\alpha -1}> 0$ and
$\tau(t_2)=T+t_2 =\frac{\alpha(t_2-t_1)}{\alpha -1}>0$. }

Under the above choice of $T$, one can easily prove that
$\frac{\tau(t_2)^\frac{n}{2}}{V(\mathbf{g}(t_2))}=
\frac{\tau(t_1)^\frac{n}{2}}{V(\mathbf{g}(t_1))}$. Hence, there
exists a $t\in [t_1, t_2]$, such that
$0=\frac{d}{dt}\log\frac{\tau^\frac{n}{2}}{V}=\frac{n}{2\tau}+\frac{\int
Rd\mu}{V}=\frac{\int (R+\frac{n}{2\tau})d\mu}{V}\ge \inf \int(R +
\frac{n}{2\tau} + \Delta f)e^{-f}d\mu $. Equivalently, $
\mu_e(\mathbf{g}(t),\tau(t))\le 0$. Recall (\ref{No expanding
breathers theorem thm1 proof equ5}) and (\ref{No expanding breathers
theorem thm1 proof equ2}), we obtain
$0\le\mu_e(\mathbf{g}(t_1),\tau(t_1)) \le
\mu_e(\mathbf{g}(t),\tau(t))\le 0$, i.e.
$\mu_e(\mathbf{g}(t_1),\tau(t_1)) = 0$. Hence all the inequalities
in (\ref{No expanding breathers theorem thm1 proof equ2}) must be
equalities and the first variation of $\mathcal{W}_e$ must vanish.
Therefore, the Ricci breather must be a gradient Ricci soliton,
namely,
\begin{equation}\label{No expanding breathers theorem thm1 proof equ6}
R_{ij}+\nabla_i\nabla_j f-\Dfrac{1}{2\tau}\mathbf{g}_{ij} = 0.
\end{equation}
\end{proof}

\begin{rema}
One cannot get similar results for shrinking breathers by using the
same method because of the ``wrong'' sign of $(\alpha -1)$ in the
shrinking case.
\end{rema}

\section{\textnormal{\bf Expanding Ricci breathers are Einstein }}\label{Section Expanding Ricci breathers are Einstein}
In this section, we prove the following result by using the
functional method. This result was first proved by T. Ivey in
\cite{I}, see also in \cite{BM}.

\begin{coro}\label{Expanding Ricci breathers are Einstein thm1}
Expanding Ricci breathers on compact Riemannian manifolds must be
Einstein.
\end{coro}

\begin{proof}({\bf Corollary \ref{Expanding Ricci breathers are Einstein thm1}})
The proof is very similar to the proof of Corollary \ref{No
expanding breathers theorem thm1}. We define
$\mu_{ek}(\mathbf{g},\tau) = \inf \cW_{ek}(\mathbf{g},\tau,f) =\inf
\tau^2\Dint_M \big[k(R+\frac{n}{2\tau})+\Delta f\big]e^{-f}d\mu $,
where the infimum is taken over all the smooth functions satisfying
$\int_Me^{-f}d\mu = 1$. By solving the backward heat equation in
Corollary \ref{Entropy formulae over Expanders cor2}, we get a
smooth function $f: M\times [t_1, t_2]\longrightarrow \mathbb{R}$
with initial condition $f(t_2) = f$. We define a linear function
$\tau:[t_1,t_2]\longrightarrow (0,+\infty)$ with $\tau(t_2)=T+t_2$,
where $T$ is a real number. Under the coupled system (\ref{Entropy
formulae over Expanders cor1 equ1}), by the monotonicity formula, we
obtain the following
\begin{equation}\label{Expanding Ricci breathers are Einstein thm1 proof equ2}
\begin{array}{rll}
\mu_{ek}(\mathbf{g}(t_1),\tau(t_1)) &= \inf\big|_{t = t_1}\cW_{ek}(\mathbf{g},\tau,f)\\
& \le \cW_{ek}(\mathbf{g},\tau,f)\big|_{t = t_1}\\
& \le \cW_{ek}(\mathbf{g},\tau,f)\big|_{t = t_2}\\
&=\mu_{ek}(\mathbf{g}(t_2),\tau(t_2)).
\end{array}
\end{equation}
On the other hand, it is easy to see that, if we simultaneously
scale $\tau$ and $\mathbf{g}$ by a scalar $\alpha>1$, we have the
following scaling property of $\mu_{ek}$,
\begin{equation}\label{Expanding Ricci breathers are Einstein thm1 proof equ3}
\mu_{ek}(\alpha\mathbf{g}(t),\alpha\tau) =
\alpha\mu_{ek}(\mathbf{g}(t),\tau).
\end{equation}
Combining (\ref{Expanding Ricci breathers are Einstein thm1 proof
equ2}), (\ref{Expanding Ricci breathers are Einstein thm1 proof
equ3}), and the diffeomorphic invariant property of the functionals,
we have
\begin{equation}\label{Expanding Ricci breathers are Einstein thm1 proof equ4}
\mu_{ek}(\mathbf{g}(t_1),\tau(t_1)\le
\mu_{ek}(\mathbf{g}(t_2),\tau(t_2)) =
\alpha\mu_{ek}(\mathbf{g}(t_1),\tau(t_1)).
\end{equation}
This yields $0\le (\alpha -1)\mu_{ek}(\mathbf{g}(t_1),\tau(t_1))$,
i.e.
\begin{equation}\label{Expanding Ricci breathers are Einstein thm1 proof equ5}
\mu_{ek}(\mathbf{g}(t_1),\tau(t_1))\ge 0.
\end{equation}

As in the previous section, we need to choose $T = \frac{t_2-\alpha
t_1}{\alpha -1}$ in order to simultaneously change $\mathbf{g}$ and
$\tau$. Again, one can show that
$\frac{\tau(t_2)^\frac{n}{2}}{V(\mathbf{g}(t_2))}=
\frac{\tau(t_1)^\frac{n}{2}}{V(\mathbf{g}(t_1))}$. There exists a
$t\in [t_1, t_2]$, such that
$0=\frac{d}{dt}\frac{\tau^\frac{n}{2}}{V}$. Hence, at this $t$,
$0=\frac{d}{dt}k\log\frac{\tau^\frac{n}{2}}{V}=k(\frac{n}{2\tau}+\frac{\int
Rd\mu}{V})= \frac{\int k(R+\frac{n}{2\tau})d\mu}{V}\ge \inf
\int\big[k(R + \frac{n}{2\tau}) + \Delta f\big]e^{-f}d\mu $.
Equivalently, $ \mu_{ek}(\mathbf{g}(t),\tau(t))\le 0$. Recall
(\ref{Expanding Ricci breathers are Einstein thm1 proof equ5}) and
(\ref{Expanding Ricci breathers are Einstein thm1 proof equ2}), we
obtain $0\le\mu_{ek}(\mathbf{g}(t_1), \tau(t_1)) \le
\mu_{ek}(\mathbf{g}(t),\tau(t))\le 0$, i.e.
$\mu_{ek}(\mathbf{g}(t_1),\tau(t_1)) = 0$. Hence all the
inequalities in (\ref{Expanding Ricci breathers are Einstein thm1
proof equ2}) must be equalities and the first variation of
$\cW_{ek}$ must vanish. Therefore, the Ricci breather must be a
gradient Ricci soliton, and furthermore, it must be Einstein, namely
\begin{equation}\label{Expanding Ricci breathers are Einstein thm1 proof equ6}
\begin{array}{rll}
R_{ij}+\nabla_i\nabla_j f-\frac{1}{2\tau}\mathbf{g}_{ij} &= 0\\
{\rm and}\quad R_{ij}-\frac{1}{2\tau}\mathbf{g}_{ij} &= 0.
\end{array}
\end{equation}
\end{proof}
\begin{rema}
One cannot get similar results for shrinking breathers by using the
same method because $\alpha -1< 0$ in the shrinking case.
\end{rema}

\section{Acknowledgement}
We would like to express our gratefulness to professor Shihshu
Walter Wei and professor Meijun Zhu for their constant support and
encouragement. We are indebted to professor Bennett Chow for sharing
draft of his new books and his encouragement and many inspiring
conversations . We also thank professor Xiaodong Cao and Yilong Ni
for helping us improve our first draft and valuable suggestions. We
thank MSRI for providing financial support during the visit.

Finally, we would like to thank the referee for helpful suggestions
which simplify the original writing greatly.


\begin{thebibliography}{99}
\bibitem{C}Cao, Xiao-Dong. {\em Eigenvalues of ($-\Delta + \dfrac{R}{2}$) on manifolds with
nonnegative curvature operator.} Mathematische Annalen (to appear).

\bibitem{CL}Chang, Shu-Cheng; Lu, Peng. {\em Evolution of Yamabe constant under Ricci
flow.} arxiv:math.DG/0604546.

\bibitem{BD}Chow, Bennett; Knopf, Dan. {\em The Ricci flow: An introduction}.
Mathematical Surveys and Monographs, AMS, Providence, RI, 2004.

\bibitem{BM}Chow, Bennett; Chu, Sun-Chin; Glickenstein, David;
Guenther, Christine; Isenberg, Jim; Ivey, Thomas; Knopf, Dan; Lu,
Peng; Luo, Feng; Ni, Lei. {\em The Ricci Flow: Techniques and
Applications. Volume 2 - Part I Geometric Aspects.} Preprint of a
book.

\bibitem{BNL}Chow, Bennett; Lu, Peng; Ni, Lei. {\em Hamilton's Ricci flow}.
Preprint of a book.

\bibitem{FIN}Feldman, Mikhail; Ilmanen, Tom; Ni, Lei. {\em
Entropy and reduced distance for Ricci expanders}. (English. English
summary) J. Geom. Anal. 15 (2005), no. 1, 49-62.

\bibitem{GT}Gilbarg, David; Trudinger, Neil S. {\em Elliptic Partial
Differential Equations of Second Order}. (Second edition.)
Springer-Verlag, 1983.

\bibitem{H1}Hamilton, Richard S. {\em Three-manifolds with positive
Ricci curvature}. J. Differential Geom., 17(2): 255-306, 1982.

\bibitem{H}Hamilton, Richard S. {\em Four-manifolds with positive curvature operator}.
J. Differential Geom., 24(2): 153-178, 1986.

\bibitem{Hu}Hsu, Shu-Yu. {\em A simple proof on the non-existence of shrinking breathers for the
Ricci flow}. Calc. Var. Partial Differential Equations 27 (2006),
no. 1, 59-73.

\bibitem{I}Ivey, Thomas. {\em Ricci solitons on compact three-manifolds}.
Differential Geom. Appl. 3 (1993), no. 4, 301--307.

\bibitem{KZ}Kuang, Shilong; Zhang, Qi S. {\em A gradient estimate for all positive
solutions of the conjugate heat equation under Ricci flow}. Personal
communications.

\bibitem{M}Ma, Li. {\em Eigenvalue monotonicity for the Ricci-Hamilton flow}.
Ann. Global Anal. Geom. 29 (2006), no. 3, 287--292.

\bibitem{N}Ni, Lei. {\em The entropy formula for linear heat equation}. J. Geom. Anal. 14
(2004) 87¨C100; Addendum. J. Geom. Anal. 14 (2004) 369-374.

\bibitem{P}Perelman, Grisha. {\em The entropy formula for the Ricci flow
and its geometric applications}. arXiv:math.DG/0211159.

\bibitem{T}Topping, Peter. {\em Lectures on the Ricci
Flow},\\ http://www.maths.warwick.ac.uk/~topping/RFnotes.html
\end{thebibliography}
%

\end{document}